\documentclass[11pt,reqno]{amsart}
\usepackage{amsfonts}
\usepackage{amsmath}
\usepackage{amssymb}
\usepackage{verbatim}        % comment the following line to get rid
\newtheorem{theorem}{Theorem}
\newtheorem{prop}{Proposition}

\newtheorem{lemma}[prop]{Lemma}

\newcommand{\ifs}{\mathcal{I}}
\newcommand{\mr}{\mathcal{R}}
\newcommand{\ol}{\overline}

\newcommand{\mg}{\mathcal{G}}
\newcommand{\e}{\varepsilon}
\newcommand{\leb}{\mathcal{L}}

\DeclareMathOperator{\diam}{diam}
\DeclareMathOperator{\dist}{dist}
\DeclareMathOperator{\id}{Id}

\title{Moreira's Theorem on the arithmetic sum of dynamically defined Cantor sets}
\author{Pablo Shmerkin}
\address{School of Mathematics (Alan Turing Building)\\
         University of Manchester\\
         Manchester\\ 
         M13 9PL\\
         United Kingdom
         }
\email{Pablo.Shmerkin@manchester.ac.uk}

\begin{document}

\begin{abstract}
We present a complete proof of a theorem of C.G. Moreira. Under mild
checkable conditions, the theorem asserts that the Hausdorff
dimension of the arithmetic sum of two dynamically defined Cantor
subsets of the real line, equals either the sum of the dimensions
or $1$, whichever is smaller.
\end{abstract}

\maketitle

\section{Introduction}

A classical problem in geometric measure theory is to compute or
estimate the fractal dimension of the arithmetic sum $K^1 + K^2$
in terms of the dimensions of $K^1, K^2\subset\mathbb{R}^n$. The
sumset $K^1+K^2$ is, up to affine equivalence, the orthogonal
projection of the product set $K^1\times K^2$ onto the line $\{
(t,t)\}_{t\in\mathbb{R}}$, so this problem is related to questions
on orthogonal projections. See \cite{PeresShmerkin} for further
discussion on this connection and the history of the problem. 

Since orthogonal projections are
Lipschitz maps and do not increase dimension, it follows that
\[
\dim_H(K^1+K^2) \le \dim_H(K^1\times K^2) \le \ol{\dim}_B(K^1) +
\dim_H(K^2),
\]
where $\dim_H$ denotes Hausdorff dimension and $\ol{\dim}_B$
denotes upper box (Minkowski) dimension; see \cite{mattila} for
the right-hand side inequality. If $\dim_H(K^1)=\dim_B(K^1)$
(which is the case if, for example, $K^1$ is the attractor of  a
self-conformal iterated function system), we obtain the inequality
\begin{equation} \label{eq:sumineq}
\dim_H(K^1+K^2) \le \min(\dim_H(K^1) + \dim_H(K^2),n),
\end{equation}
where $n$ is the dimension of the ambient space. Obtaining lower
bounds is much harder, and it is easy to construct examples where
(\ref{eq:sumineq}) fails. However, there is a heuristic principle
which says that ``generically'' (\ref{eq:sumineq}) should hold as
equality. In some cases this has been accomplished in a
measure-theoretical sense, see for example
\cite{peres-solomyak-cantor}. However, from those results one
cannot tell whether equality in (\ref{eq:sumineq}) holds for a
{\em specific} pair $K^1$ and $K^2$.

An \textbf{iterated function system} (or i.f.s. for short) is a
finite family $\{ f_1,\ldots, f_m\}$ of self-maps of
$\mathbb{R}^n$ (or a more general complete metric space, but here
we will only consider iterated function systems on the real line),
such that each map $f_i$ is Lipschitz with Lipschitz constant
strictly less than $1$; in other words, such that
\[
|f_i(x) - f_i(y)| \le L |x-y| \quad \textrm{ for all }
x,y\in\mathbb{R}^n, \, i\in\{1,\ldots,m\},
\]
for some constant $L<1$. If all the maps $f_i$ are $C^\alpha$ for
some $\alpha\ge 1$ we will say that the i.f.s. is $C^\alpha$. For
fixed $\alpha$ and $m$, the family of all $C^\alpha$ iterated
function systems with $m$ maps inherits a natural topology from
$C^\alpha\times \cdots\times C^\alpha$, where the product is of
course $m$-fold.

Given an i.f.s. $\ifs=\{f_1,\ldots, f_m\}$, the \textbf{attractor}
$K=K(\ifs)$ is the only nonempty compact subset of $\mathbb{R}^n$
such that
\[
K = \bigcup_{i=1}^m f_i(K).
\]
See e.g. \cite{falconer2} for more background on iterated function
systems, including the existence and uniqueness of attractors.

A \textbf{regular Cantor set} $K\subset\mathbb{R}$ is the
attractor of a $C^2$ iterated function system $\{f_1,\ldots,
f_m\}$ such that the sets $f_i(I)$ are pairwise disjoint, where
$I$ is the convex hull of $I$, and moreover $f_i: I\rightarrow
f_i(I)$ is a diffeomorphism. In the dynamics literature regular
Cantor sets are usually defined as repellers of smooth expanding
maps. The existence of Markov partitions allows to realize such
repellers as attractors of iterated function systems of a more
general kind (i.e. ``graph directed'' ones). In this paper we
concentrate on the most basic kind of attractors, but this is just
a matter of notational simplicity; both the result and the proof
extend in a straightforward way to more general repellers.

Moreira and Yoccoz \cite{moreirayoccoz} proved a deep result about
the arithmetic sum of regular Cantor sets $K^1, K^2$ when
$\dim(K^1)+\dim(K^2)>1$. They prove that generically (in a
topological sense with respect to the $C^2$ topology) such
sumsets contain intervals, settling a conjecture of Jacob Palis
\cite{palis-conjecture}. The results in \cite{moreirayoccoz} have
important consequences on the study of homoclinic bifurcations.

In a different direction,  Moreira \cite{moreira-hungarica} studied
some problems in diophantine approximation which also involve sums
of Cantor sets. As part of the solution to those problems, Moreira
states the following result:
\begin{theorem} \label{th:gugu}
Let $\{f^i_1,\ldots, f^i_{m_i}\}$, $i=1,2$, be a pair of $C^2$
iterated function systems on $\mathbb{R}$, and let $K^1, K^2$ be
the attractors. Suppose that the families $\{f^i_j\}_{j=1}^{m_i}$,
$i\in\{1,2\}$ are pairwise disjoint. Assume that the following
properties hold:
\begin{enumerate}
\item There are $1\le i < j \le m_1$ and
$x_0$ in $K^1$ such that
\[
\left(f^1_i \circ (f^1_j)^{-1}\right)''(x_0) \neq 0.
\]
\item There exist $1\le l_i \le m_i$, $i=1,2$, such that if $y_i$ is the fixed point of $f^i_{l_i}$ then
\[
\frac{\log |(f^1_{l_1})'(y_1)|}{\log |(f^2_{l_2})'(y_2)|} \notin
\mathbb{Q}.
\]
\end{enumerate}
Then
\[
\dim_H(K^1+K^2) = \min(\dim_H(K^1)+\dim_H(K^2),1).
\]
\end{theorem}
The hypotheses in this theorem are generic (the first one is
robust and dense in the $C^1$ topology, while the second holds for
almost every parameter in generic parametrized families).
Moreover, the hypotheses are explicit and can be checked in
specific examples.

Unfortunately, the proof of Theorem \ref{th:gugu} which appeared
in \cite{moreira-hungarica} has some errors. Even though the basic idea is correct, it is far from trivial to recover a complete proof for it, and even the basic ideas may be somewhat obscure for those not familiar with the techniques in
\cite{moreirayoccoz}. C.G. Moreira explained to us the main corrections
needed; based on this we were able to reconstruct a complete proof
of Theorem \ref{th:gugu}. The purpose of this note is to write down this proof in detail. One motivation for doing this is that we believe that some of the ideas contained in the proof may find application in other problems in geometric measure theory or dyamics,
where current methods only yield almost everywhere or random results.

Several developments took place after a first version of this
paper was completed. Moreira informed us that he can now prove
Theorem \ref{th:gugu} without assuming hypothesis (2). Using some
of the ideas presented in this paper, but also substantial new
ones, Y. Peres and the author \cite{PeresShmerkin} proved a
version of Theorem \ref{th:gugu} when the $f^i_j$ are all linear
maps; this includes classical examples such like central Cantor
sets. We also prove that hypothesis (2) is necessary in this case.
Ero\v{g}lu \cite{eroglu-sums} investigated the Hausdorff measure
of sumsets in the critical dimension. In particular, he proves
that in many cases it is zero. Thus we now have a rather complete
picture of the size of the arithmetic sum of dynamically-defined
Cantor sets in the line, at least in the case where the sum of their dimensions does not exceed one.

\section{Notation}

Let $\ifs=\{f_1,\ldots,f_m\}$ be a $C^2$ i.f.s. on $\mathbb{R}$
(with $\e>0$), such that the basic pieces $f_i(K)$ are pairwise
disjoint, where $K=K(\ifs)$ is the attractor. We say that $\ifs$
is \textbf{normalized} if the convex hull of $K$ is the unit
interval $I=[0,1]$.

Orientation-preserving (surjective) diffeomorphisms of the unit interval will
play an important role. The set of all such diffeomorphisms of
class $C^1$, endowed with the $C^1$ topology, will be denoted by
$\mg$. An alternative way of thinking of $\mg$ is as the space of all
diffeomorphic embeddings of the unit interval into $\mathbb{R}$,
divided by the action of the affine group by left composition (the
equivalence of both definitions is given by the choice of a
representative in a canonical way).

We will use the following form of the $C^1$ norm:
\[
\| f\|_{C^1} = \max\{ \|f\|_{L^\infty},\|f'\|_{L^\infty} \}.
\].
In addition, we let
\[
\mg(\delta) = \{ g\in\mg: \|g-\id\|_{C^1}<\delta \}.
\]

We record the following immediate lemma for later
reference:
\begin{lemma} \label{lemma:norm}
Let $f:J\rightarrow \mathbb{R}$ be a diffeomorphism, where $J$ is
a closed subinterval of $I$. Then $\|f-\id\|_{C^1} =
\|f'-1\|_{L^\infty}$.
\end{lemma}
\textit{Proof}. Obviously $\|f-\id\|_{C^1}\ge
\|f'-1\|_{L^\infty}$. The other inequality also follows since for
$x\in I$ we have
\[
|f(x)-x| \le \int_J |f'(x)-1| dx \le \|f'-1\|_{L^\infty}.
\]
\qed

Let $\ifs$ be any regular normalized i.f.s. with attractor $K$.
The symbolic space is $\Sigma=\{1,\ldots,m\}^\mathbb{N}$, where
$m$ is the number of maps in the i.f.s. The set of all finite
words with symbols in $\{1,\ldots,m\}$ will be denoted by
$\Sigma^*$.

If $u=(u_1,\ldots, u_j)\in\Sigma^*$, we will write $f_u =
f_{u_1}\circ\cdots f_{u_j}$. The reverse word $(u_j,\ldots, u_1)$
will be denoted by $u^\star$. We will also let $T_u$ be the unique
affine map such that $T_u\circ f_{u}\in\mg$.

Given $\omega\in\Sigma$, let $\omega|k$ be the restriction of
$\omega$ to the first $k$ coordinates. For a given
$\omega\in\Sigma$, consider the sequence
$T_{(\omega|k)^\star}\circ f_{(\omega|k)^\star}\in\mg$. Sullivan
\cite{sullivan} proved that this sequence converges in $C^1$,
uniformly in $\omega$ (in fact the convergence is in any
smoothness class to which the $f_i$ belong, but for us $C^1$
suffices). The limiting diffeomorphism will be denoted by
$L_\omega$, and the image set $L_{\omega}(K)$ will be called a
\textbf{limit geometry} of $K$. Limit geometries are also regular,
normalized Cantor sets with the same dimension as $K$ (indeed,
$L_\omega K$ is the attractor of $\{ L_\omega f_i
L_\omega^{-1}\}_{i=1}^m$). Moreover, since the convergence in
Sullivan's Theorem is uniform and $\Sigma$ is compact, the family
$\{ L_\omega: \omega\in\Sigma\}$ is also compact.

We will consider pairs of iterated function systems $\ifs^1,
\ifs^2$ on $\mathbb{R}$, and the product attractor $\Lambda =
K^1\times K^2=K(\ifs^1)\times K(\ifs^2)$. Throughout the paper we
will distinguish the i.f.s. we are referring to by the use of a
superscript. For example, $\Sigma^1, \Sigma^2$ will denote the
symbol spaces corresponding to the i.f.s. $\ifs^1,\ifs^2$
respectively. It should be clear from the context whether a
superscript is used in this fashion, or to denote a power
operation.

We will always denote $d^i=\dim(K^i), i=1,2$, and $d=d^1+d^2$. The
images $f^i_{u}(K^i)$, where $u\in \Sigma^i$, will be referred to
as \textbf{cylinder sets}, and denoted by $K^i(u)$. We will deal
with the convex hull of cylinder sets rather often; the convex
hull of $K^i(u)$ will be denoted by $I^i(u)$.

Let $\rho>0$ be a small number. The \textbf{$\rho$-decomposition} of $\Lambda$, denoted by
$\Lambda(\rho)$, is the collection of all pairs of words
$(u_1,u_2)$ such that
\[
|I^i(u_i)| = \diam(K^i(u_i)) > \rho  \quad (i=1,2),
\]
but these inequalities fail for any words containing $u_1, u_2$ as
proper initial subwords. For $\phi=(\phi_1,\phi_2)\in\mg\times\mg$
we will also abbreviate
\[
Q^\phi(u_1, u_2) = \phi_1(I^1(u_1))\times \phi_2(I^2(u_2)).
\]
When $\phi=\id\times\id$ (where $\id$ is the identity map) we will
simply write $Q(u_1, u_2)$.

Many calculations will depend on a previously fixed constant $A$.
Given two positive quantities $x,y$, by $x\lesssim y$ we will mean
$x < C y$ for some constant $C$ which depends continuously on $A$,
$\ifs^1$, and $\ifs^2$. We define $x \gtrsim y, x \approx y$
analogously.

When we need to refer to constants explicitly we will denote them
by $c$ or $C$; their value can be different at each line, and they
always depend continuously on $A$, $\ifs^1$, and $\ifs^2$.

Let $\Pi_\lambda:\mathbb{R}^2\rightarrow \mathbb{R}$ be the
projection-type mapping $\Pi(x,y)=x+\lambda y$. Let $\mr$ be a
subset of the $\rho$-decomposition $\Lambda(\rho)$. We will say
that $\mr$ is \textbf{$(\eta,\lambda,\phi)$-faithful} if it
contains a subfamily $\mr'$, with $\#\mr'> \rho^{\eta-d}$, and
such that
\[
\{\Pi_\lambda(Q^\phi(u_1,u_2)): (u_1,u_2)\in\mr' \}
\]
is a pairwise disjoint family. The following lemma, although very
simple, will play a crucial role in the proof:

\begin{lemma} \label{lemma:rhoperturb} Fix $C_0, \eta>0$. Then for
all sufficiently small $\rho$ (depending on $C_0$ and $\eta$) the
following holds: if a family $\mr\subset\Lambda(\rho)$ is
$(\eta,\lambda,\phi)$-faithful, it is also
$(2\eta,\widetilde{\lambda},\widetilde{\phi})$-faithful for all
$\widetilde{\lambda},\widetilde{\phi}$ such that
\[
\left|\widetilde{\lambda}-\lambda\right| \le C_0 \rho,\quad
\left\|\widetilde{\phi}_i-\phi_i\right\|_{C^1} \le C_0
\rho\,(i=1,2).
\]
\end{lemma}
\textit{Proof}. Note that there exists $C=C(C_0,\eta)>0$ such that
\[
\Pi_{\widetilde{\lambda}}\left(Q^{\widetilde{\phi}}(u_1,u_2)\right)\subset
C\cdot \Pi_{\lambda}(Q^\phi(u_1,u_2)),
\]
where $C\cdot J$ denotes the interval with the same center as $J$
and length $C|J|$. Therefore if $\mr_1$ is the family arising from
the definition of $(\eta,\lambda,\phi)$-faithful, there is a
subset $\mr_2$ of $\mr_1$ of cardinality at least $|\mr_1|/(2 C)$
such that
\[
\{ C \cdot \Pi_{\lambda}(Q(u_1,u_2)): (u_1,u_2)\in\mr_2 \}
\]
is a disjoint family. Taking $\rho$ small enough so that
$\rho^\eta < (2C)^{-1}$ yields the lemma. \qed

Finally, we define the renormalization operators. This is a family
of operators $\{R_{u_1, u_2}:u_i\in(\Sigma^*)^i\}$, defined as
\[
R_{u_1,u_2}(\omega_1, \omega_2, s) = \left(u_1 \omega_2, u_2
\omega_2,
 \left|I^1(u_1)\right|^{-1} \left|I^2(u_2)\right|s\right).
\]

To understand the action of these operators, define
\[
\Lambda_{\omega_1,\omega_2} = L_{\omega_1}\left(K^1\right) \times
L_{\omega_2}\left(K^2\right).
\]
We will also refer to $\Lambda_{\omega_1,\omega_2}$ as a limit
geometry of $\Lambda$. Small cylinders are very close, after
rescaling, to a limit geometry, and since we will deal with robust
properties (i.e. properties which are invariant under small
perturbations of the parameters) we will able to draw conclusions
about limit geometries from its finite approximations.

Moreover, a cylinder of a limit geometry is also close to a
corresponding cylinder in the original attractor. In this sense,
the action of the renormalization operator $R_{u_1,u_2}$ is to
``zoom into'' the $(u_1,u_2)$-cylinder of the given limit geometry
(or approximating cylinder); the transformation of $s$ simply
takes into account the normalization (rescaling back to the unit
square) of the cylinder.

\section{Auxiliary results}

In this section we collect a number of basic results that we will
use in the course of the proof of the main theorem.

\begin{lemma} \label{lemma:linear}
Let $\ifs=\{f_i\}_{i=1}^m$ be a regular i.f.s. with attractor $K$.
Then the following holds:
\begin{enumerate}
\item[(i)]
If $f_1$ is linear then for all finite words $u$ the limit
geometry $L_{u 1^\infty}(K)$ is affinely equivalent to the
cylinder $f_{u^\star}(K)$. More precisely, we have the identity
\[
L_{u 1^\infty} = T_{u^\star} f_{u^\star}.
\]
\item[(ii)] The map $g_j = L_{j^\infty}f_j L_{j^\infty}^{-1}$ is linear, and its eigenvalue is equal to the eigenvalue of
$f_j$.
\item[(iii)] Suppose that $f_1$ is linear, and fix a finite word $u$ with symbols in $\{1,\ldots,m\}$.
Consider a new i.f.s $\{ h_1,\ldots, h_m\}$, where
\[
h_i = L_{u 1^\infty} f_i L_{u 1^\infty}^{-1}.
\]
Then for all words $v$,
\[
|h_{v^\star}(I)| = \frac{|I((vu)^\star)|}{|I(u^\star)|}.
\]
\end{enumerate}
\end{lemma}
\textit{Proof}.
\begin{enumerate}
\item[(i)] This is obvious when thinking of $\mg$ as a quotient space.
\item[(ii)] Again using the quotient space interpretation, it is clear
that the class of $L_{j^\infty} f_j$ is the same as the class of
$L_{j^\infty} $, and therefore the class of $L_{j^\infty}f_j
L_{j^\infty}^{-1}$ is the affine group. The invariance of
eigenvalues under conjugacies is a general fact.
\item[(iii)]
Keeping in mind that $T_z$ is linear and using (i) we have:
\begin{eqnarray*}
|h_{v^\star}(I)| & =  & | T_{u^\star} f_{u^\star}
f_{v^\star} (T_{u^\star} f_{u^\star})^{-1} I| \\
& = & | T_{u^\star} f_{u^\star} f_{v^\star} I | =  | T_{u^\star}
f_{(vu)^\star} I| \\
& = & | T_{u^\star} I | | T_{(vu)^\star} I |^{-1}
=|I(u^\star)|^{-1} |I((vu)^\star)|.
\end{eqnarray*}\qed
\end{enumerate}

The previous lemma will allow us to assume that the maps in $K^1,
K^2$ for which the incommensurability holds are actually linear,
and this in turn will imply that limit geometries are cylinder
sets.

We will need the well-known
\textbf{bounded distortion principle} (see e.g. \cite[Proposition 4.2]{falconer3} :
\begin{lemma}
Let $\ifs$ be a $C^{1+\e}$ i.f.s. for some $\e>0$, and let $K$ be
the attractor. Then there is $L_1>0$ such that
\[
L_1^{-1} < \frac{|f'_u(x)|}{|f'_u(y)|} < L_1,
\]
for all $u\in\Sigma$ and all $x,y\in K$.
\end{lemma}
The following proposition is a key geometrical result. It is a
kind of discrete Marstrand theorem on projections in a particular
setting (see \cite[Chapter 10]{mattila} for general projection
theorems). Intersection numbers are used in a similar fashion in
the work of M. Rams, see \cite{rams-pacific} and references
therein. In \cite{PeresShmerkin} a general discrete projection
theorem is proved, but the result we need here does not follow
directly from it, so a full proof is given for the convenience of
the reader.

\begin{prop} \label{prop:fixedscale}
Fix a large constant $A$ and a small constant $\eta>0$. Let $K^1,
K^2$ be attractors of regular normalized i.f.s. $\ifs^1, \ifs^2$
respectively, $d^i = \dim_H(K^i)$, $\Lambda=K^1\times K^2$ and
$d=d^1+d^2=\dim_H(\Lambda)$. Assume that $d<1$.

Then there is a number $\rho_0>0$, which depends continuously on
$\eta$, $A$, $\ifs^1$ and $\ifs^2$, such that for all
$0<\rho<\rho_0$ there exists a set $J\subset [-A,A]$ with the
following properties:
\begin{enumerate}
\item
$ \leb([-A,A]\backslash J) < \rho^\eta$, where $\leb$ denotes
one-dimensional Lebesgue measure.
\item
If $\lambda\in J$ and $\mr$ is any subset of the
$\rho$-decomposition such that $\#\mathcal{R}
> \rho^{\eta-d}$, then $\mr$ is $(4\eta,\lambda,\id)$-faithful.
\end{enumerate}

\end{prop}

\textit{Proof}. Let
\[
N(\lambda) = \#\left\{ (u,v)\in\Lambda(\rho)\times\Lambda(\rho):
\Pi_\lambda(Q(u))\cap \Pi_\lambda(Q(v)) \neq \varnothing\right\}.
\]
Note that in the above $u$ and $v$ are \textit{pairs} of words.
Given $u,v\in\Lambda(\rho)$ let also
\[
E(u,v) = \{ \lambda : \Pi_\lambda(Q(u))\cap \Pi_\lambda(Q(u)) \neq
\varnothing \}.
\]
Observe that if $\lambda\in E(u,v)$ then there is a line with
slope $\lambda$ intersecting both $Q(u)$ and $Q(v)$. Therefore we
have the estimate
\begin{equation} \label{eq:estcommondirs}
\leb(E(u,v)) \lesssim \rho/\dist(Q(u),Q(v)).
\end{equation}
As a consequence of the bounded distortion principle, the
following holds: given $u\in\Lambda(\rho)$ and $\e>\rho$,
\begin{equation} \label{eq:estnumberrects}
\# \{ v\in\Lambda(\rho) : \dist(Q(u),Q(v)) < \varepsilon \}
\lesssim
 (\rho/\varepsilon)^{-d}.
\end{equation}
The constant implied by the $\lesssim$ notation depends
continuously on $\ifs^1, \ifs^2$. In particular, $\#\Lambda(\rho)
\lesssim \rho^{-d}$, and $\dist(Q(u),Q(v)) \gtrsim \rho$ for any
two different $u, v$ in the $\rho$-decomposition.

\textbf{Claim}.
\[
\int_{-A}^A N(\lambda) d\lambda \lesssim \rho^{-d}.
\]
\textit{Proof of Claim}. Given $u,v\in\Lambda(\rho)$ let
$\textrm{d}(u,v) = \dist(Q(u),Q(v))$. We estimate:
\begin{eqnarray*}
\int_{-A}^A N(\lambda) d\lambda & = & \sum_{u\in\Lambda(\rho)} \sum_{v\in\Lambda(\rho)} \int_{-A}^A \mathbf{1}_{ \{\Pi_\lambda(Q(u))\cap \Pi_\lambda(Q(v))\neq \varnothing \}} d\lambda\\
& =  & \sum_{u \in \Lambda(\rho)} \sum_{i=0}^{-\log_2(C \rho)}
\sum_{1 < 2^ \textrm{d}(u,v)\le 2} \leb(E(u,v)) \\
& \lesssim & \rho \sum_{u\in\Lambda(\rho)} \sum_{i=0}^{-\log_2(C
\rho)} 2^i \#\{v\in\Lambda(\rho): \textrm{d}(u,v) < 2^{1-i}\}\\
& \lesssim &  \rho \sum_{u\in\Lambda(\rho)}
\sum_{i=0}^{-\log_2(C \rho)} 2^i (\rho 2^{i-1})^{-d} \\
& \lesssim & \rho^{1-d} \#\Lambda(\rho)
2^{-(1-d)\log_2(C\rho)} \\
& \approx &  \rho^{-d},
\end{eqnarray*}
where we used (\ref{eq:estcommondirs}) in the third line and
(\ref{eq:estnumberrects}) in the fourth line. This proves the
claim.
Let $J$ be defined as
\[
J = \{ \lambda\in [-A,A]: N(\lambda) < \rho^{-2\eta-d}\}.
\]
We will show that $J$ has the desired properties. Firstly, by the
claim and Chebychev's inequality,
\[
\leb([-A,A]\backslash J) \lesssim \rho^{2\eta}\,\Longrightarrow \,
\leb([-A,A]\backslash J) < \rho^{\eta},
\]
if $\rho$ is small enough. Now let $\mr$ be a subset of
$\Lambda(\rho)$ such that $\#\mr
> \rho^{\eta-d}$, and define
\[
N_1(\lambda) = \#\{ (u,v)\in\mr \times \mr : \Pi_\lambda(Q(u))\cap
\Pi_\lambda(Q(v))\neq \varnothing \},
\]
and define $J_1$ analogously using $N_1$ instead of $N$. Clearly
$N_1(\lambda) \le N(\lambda)$ for all $\lambda$, whence $J\subset
J_1$. Therefore it is enough to prove that condition (2) applied
to $\mr$ in the proposition holds for any fixed $\lambda\in J_1$.

Note that if $\lambda\in [-A,A]$ then $\Pi_\lambda(\Lambda)\subset
[-A-1,A+1]$. Let us divide $[-A-1,A+1]$ into intervals $I_j$ of
length slightly less than $\rho$, $j=1,\ldots, \lceil
(2A+2)\rho^{-1} \rceil$. Write $m_j$ for the number of rectangles
$Q=Q(u_1,u_2)$, where $(u_1,u_2)\in\mr$, such that the center of
$I_j$ belongs to $\Pi_\lambda(Q)$. Also let
\[
\mathcal{J} = \{ 1\le j \le \lceil (2A+2)\rho^{-1} \rceil : m_j >
0 \},
\]
and observe that it is enough to show that $\#\mathcal{J} \gtrsim
\rho^{4\eta-d}$ (provided this holds, for each $j\in\mathcal{J}$
we pick $(u_1,u_2)\in\mr$ such that the center of $I_j$ belongs to
$\Pi_\lambda(Q(u_1,u_2))$; by construction this is a family with a
bounded covering number, so we can pick an appropriate disjoint
subfamily $\mr'$ with comparable cardinality).

Note that each $\Pi_\lambda(Q)$ contains the center of a uniformly
bounded number of $I_j$, and therefore
\[
\sum_{j\in\mathcal{J}} m_j \gtrsim \#\mr \ge \rho^{\eta-d}.
\]
Using this we estimate, for sufficiently small $\rho$,
\begin{eqnarray*}
\rho^{-2\eta-d} & > & N_1(\lambda) \ge  \sum_{j\in\mathcal{J}} m_j^2 \\
& \ge & (\#\mathcal{J})^{-1} \left(\sum_{j\in\mathcal{J}}
m_j\right)^2 \ge (\#\mathcal{J})^{-1} \rho^{2\eta-2d}.
\end{eqnarray*}
This concludes the proof of the proposition. \qed

We remark that because of the compactness of the set of limit
geometries, the number $\rho_0$ given by the proposition can be
chosen uniformly for all limit geometries
$\Lambda_{\omega_1,\omega_2}$.

\section{Proof of the main theorem}

\subsection{Sketch of proof}

We begin by sketching the proof; full details follow below. The
bulk of the proof consists in showing that given $\varepsilon>0$,
the inequality
\[
\dim_H(\Pi_\lambda(\Lambda)) > d - \varepsilon
\]
holds for $\lambda$ in some {\em open} set. Moreover, to begin
with we can assume that $\ifs^1$ and $\ifs^2$ both contain a
linear map. From here one can deduce, using incommensurability,
that the same inequality holds for {\em all} $\lambda\neq
0,\infty$, and then pass to the general case by approximating
cylinders by limit geometries and using Lemma \ref{lemma:linear}.

We fix a small $\eta$ and for each $\rho>0$ apply Proposition
\ref{prop:fixedscale} to obtain many (more precisely, all up to a
small exponential correction) rectangles in the
$\rho$-decomposition with disjoint projections. For a fixed
$\rho$, this construction is robust in $\lambda$ (perturbing the
$\lambda$ slightly the rectangles will have projections at a
distance of at least, say, $\rho/2$).

The goal is to carry this construction inductively in each of
those rectangles, but \textit{a priori} there is a big obstacle:
the set of parameters $\lambda$ given by Proposition
\ref{prop:fixedscale}, even though of almost full measure, can
vary for each rectangle, so we need some device to make sure that
the parameters are {\em recurrent}; i.e. we can take the {\em
same} set $J$ for all rectangles, perhaps at the price of reducing
the number of rectangles we are working with slightly. Such
recurrence result was obtained in \cite{moreirayoccoz} and is one
of the main technical tools in the proof of Theorem \ref{th:gugu}.

This inductive construction yields, for each $\lambda$ in some
open set, a Moran construction whose limit set is contained in
$\Pi_\lambda(\Lambda)$ and has dimension at least $d-\varepsilon$,
provided $\eta$ and then $\rho$ were taken sufficiently small.
This finishes the sketch of the proof.

\subsection{The scale recurrence lemma}

For the convenience of the reader we state the key scale
recurrence lemma (sometimes called the scale selection lemma). For
the proof, the reader is referred to \cite{moreirayoccoz}.

Let us say that a regular i.f.s. $\ifs$ is \textbf{essentially
nonlinear} if it verifies condition (1) in Theorem \ref{th:gugu};
in other words, if there exist $i<j, x_0\in K$ such that
\[
\left( f_i \circ f_j^{-1} \right)''(x_0) \neq 0.
\]

\begin{theorem} \label{th:srl} Given regular Cantor sets
$K^1, K^2$, such that at least one of $K^1, K^2$ is essentially
nonlinear, there exists a large constant $A$ such that, setting
$a=A^{-1}$, the following holds:

Let $\rho$ be sufficiently small. Suppose that for each
$\omega_1\in\Sigma^1,\omega_2\in\Sigma^2$ some measurable set
$J_{\omega_1,\omega_2}$ is given such that
\[
\leb(I_A\backslash J_{\omega_1,\omega_2}) < a,
\]
where $I_A = [-A,-1/A]\cup [1/A,A]$. Then there exists another
family 
\[
\{ F_{\omega_1,\omega_2}\}_{\omega_1\in \Sigma^1,\omega_2\in\Sigma^2},
\]
verifying the following properties:
\begin{enumerate}
\item $F_{\omega_1,\omega_2}$ is contained in the
$(A\rho)$-neighborhood of $J_{\omega_1,\omega_2}$.
\item For every $s\in F_{\omega_1,\omega_2}$ there are at least $a
\rho^{-d}$ elements of the $\rho$-decomposition of
$\Lambda_{\omega_1,\omega_2}$ such that if $(u_1, u_2)$ is one
such element and
\[R_{u_1,u_2}(\omega_1,\omega_2, s) = (u_1\omega_1,u_2\omega_2, s'),
\]
then $(s'-\rho,s'+\rho)\in F_{u_1\omega_1,u_2\omega_2}$.
\end{enumerate}
\end{theorem}

\subsection{The core of the proof} We now start the proof of Theorem \ref{th:gugu}. We start by proving a weaker result; Theorem \ref{th:gugu} will be obtained later as a corollary.

\begin{prop} \label{prop:gugu2}
Let $K^1, K^2$ be regular Cantor sets of dimension $d^1, d^2$,
such that $d^1+d^2<1$ and $K^1$ is essentially nonlinear. Assume
also that the maps $f^i_1$, $i=1,2$, are linear.

Then for all $\e>0$ there exist nonempty open sets $U^+\subset\mathbb{R}^+,
U^-\subset\mathbb{R}^-$ and $\delta>0$ such that
\begin{equation} \label{eq:dimineq}
\dim_H(\phi_1(K^1)+\lambda\phi_2(K^2)) > d^1+d^2-\e
\end{equation}
for all $\lambda\in U^+\cup U^-$ and all
$\phi_1,\phi_2\in\mg(\delta)$.
\end{prop}

\textit{Proof}. Let $A$ be the constant given by the scale
recurrence lemma, and write $a=1/A$. Fix a small $\eta>0$, and
then a very small $\rho>0$ so that $\rho^\eta<a$, Proposition
\ref{prop:fixedscale} works for this $\rho$ for all limit
geometries $\Lambda_{\omega_1,\omega_2}$ and the Scale Recurrence Lemma is satisfied. In the course of the proof we will need $\rho$ to satisfy additional conditions; it will be clear that all can be satisfied by starting with a sufficiently small $\rho$.

For each pair $(\omega_1,\omega_2)\in\Sigma_1\times\Sigma_2$ let
$J_{\omega_1,\omega_2}$ be the set given by Proposition
\ref{prop:fixedscale} applied to the limit geometry
$\Lambda_{\omega_1,\omega_2}$.

We apply Theorem \ref{th:srl} to obtain a new family
$\{F_{\omega_1,\omega_2}\}$ with the conditions prescribed in the
scale recurrence lemma. Clearly if $\rho$ is small then
$F_{\omega_1,\omega_2}$ contains both positive and negative
numbers. Pick any $\lambda^+\in
F_{1^\infty,1^\infty}\cap\mathbb{R}^+$. There exists an open set
$U=U^+=(\lambda^+-c\rho,\lambda^++c\rho)$ such that the following
holds: for all $(u_1, u_2)$ arising from part (2) of the scale recurrence
lemma (applied to $F_{1^\infty, 1^\infty}$, $\lambda^+$, and
$\rho$) and all $\lambda\in U$,
\begin{equation} \label{eq:srlinside}
|I^1(u_1)|^{-1} |I^2(u_2)|\lambda \in F_{u_1 1^\infty, u_2
1^\infty}.
\end{equation}
 This follows from the fact that the quotients
$|I^2(u_2)|/ |I^1(u_1)|$ are uniformly bounded. We now fix any
$\lambda\in U$ for the rest of the proof (the construction of
$U^-$ is exactly analogous). We also fix
$\phi=(\phi_1,\phi_2)\in\mg(\rho/2)\times\mg(\rho/2)$.

We will inductively construct a tree $\mathcal{T}$, with vertices
labeled by pairs of words $(u_1, u_2)$, such that the following
holds: Let $\mathcal{T}_k$ denote the set of vertices of step $k$.
\begin{enumerate}
\item[(A)] If $(u_1, u_2)\in \mathcal{T}_{k+1}$ then $u_i = v_i z_i$ for some $(z_1, z_2)$, where $(v_1, v_2)\in
\mathcal{T}_{k}$ is the parent of $(u_1, u_2)$.
\item[(B)] If $(u_1, u_2)\in \mathcal{T}_k$ then  $|I^i(u_i)| \ge
\rho^k$, $i=1,2$. In particular,
\[
|\Pi_{\lambda}(Q^\phi(u_1,u_2))| \gtrsim \rho^k.
\]
\item[(C)]
Each vertex has $\gtrsim \rho^{9\eta-d}$ offspring.
\item[(D)] For each vertex $(u_1, u_2)\in\mathcal{T}$ the following family is pairwise disjoint:
\[
\left\{  \Pi_{\lambda}(Q^\phi(w_1,w_2)) : (u_1, u_2) \textrm{ is a
parent of } (w_1,w_2) \right\}.
\]
\end{enumerate}

Properties (A)-(D) imply that $\mathcal{T}$ induces a separated
Moran construction with cylinders $\Pi_{\lambda}(Q^\phi(u_1,
u_2))$, with limit set
\[
M = \bigcap_{k=1}^\infty \bigcup_{(u_1,u_2)\in\mathcal{T}_k}
\Pi_{\lambda}(Q^\phi(u_1, u_2)).
\]
It is clear that $M\subset \Pi_{\lambda}(\phi_1 K^1\times \phi_2 K^2)$. Moreover,
\begin{equation} \label{eq:dimbound}
\dim_H(M) \ge d-9\eta.
\end{equation}
This follows by standard methods; we sketch the proof for the convenience of the reader. We construct a probability measure $\mu$ supported on $M$ inductively as follows: suppose $\Pi_{\lambda}(Q^\phi(u_1, u_2))$ has been defined for all $(u_1,u_2)\in \mathcal{T}_k$. Then we distribute the mass of $\Pi_{\lambda}(Q^\phi(u_1, u_2))$ uniformly among all the offspring intervals $\Pi_{\lambda}(Q^\phi(v_1, v_2))$ (where $(v_1,v_2)\in\mathcal{T}_{k+1}$ ranges over the offspring of $(u_1,u_2)$). Using (B), (C) and (D), it is easy to verify that
\[
\mu(x-r,x+r) \lesssim r^{d-9\eta},
\]
for all $x\in\textrm{supp}(\mu)=M$ and all $r>0$. Thus \eqref{eq:dimbound} follows from the mass distribution principle (see \cite[Proposition 2.1]{falconer3}).

Since $\eta$ is arbitrary, it will be enough to verify properties
(A)-(D) to complete the proof.

For each $(u_1, u_2)\in \mathcal{T}_j$ we will also inductively
construct a scale $\lambda^{u_1, u_2}$ such that $ \lambda^{u_1,
u_2} \in F_{u_1^\star 1^\infty, u_2^\star 1^\infty}$ (for $j>0$).
We start by setting $\mathcal{T}_0 = \{
(\varnothing,\varnothing)\}$ (the root of the tree; here
$\varnothing$ denotes the empty word) and
$\lambda^{\varnothing,\varnothing} =\lambda$.

Now we specify the inductive construction: suppose that $(u_1,
u_2)\in \mathcal{T}_j$ for some $j$, and that $\lambda^{u_1, u_2}$
has been defined. Let $\omega_i = u_i^\star 1^\infty$, $i=1,2$,
and let us apply the scale recurrence lemma to
$\Lambda_{\omega_1,\omega_2}$ with scale $s=\lambda^{u_1, u_2}$.
We thus obtain a family of pairs of words $\mr_0^{u_1,u_2}$ given
by the scale recurrence lemma; i.e. $\#\mr_0^{u_1,u_2} > a
\rho^{-d} > \rho^{\eta-d}$, and if $(v_1^\star, v_2^\star)\in
\mr_0^{u_1,u_2}$ and we let
\[
\lambda^{v_1 u_1, v_2 u_2} = \lambda^{u_1, u_2}
\left|I^1_{\omega_1}(v_1^\star)\right|^{-1}
\left|I^2_{\omega_2}(v_2^\star)\right|,
\]
then $\lambda^{v_1 u_1, v_2 u_2} \in F_{v_1^\star u_1^\star
1^\infty, v_2^\star u_2^\star 1^\infty}$; here $I^i_{\omega_i}$
are cylinder intervals with respect to the limit geometries
$L_{\omega_i}(K^i)$. For $j=0$ this follows from
(\ref{eq:srlinside}). From Lemma \ref{lemma:linear}(iii) we get
\begin{equation} \label{eq:lambda}
\lambda^{v_1 u_1, v_2 u_2} = \frac{|I(u_1)|}{|I(u_1 v_1)|}
\frac{|I(u_2 v_2 )|}{|I(u_2)|}\lambda^{u_1, u_2}.
\end{equation}

We next use Proposition \ref{prop:fixedscale}, Lemma
\ref{lemma:rhoperturb} and the first part of the scale recurrence
lemma to obtain a subset $\mr_1^{u_1, u_2}$ of $\mr_0^{u_1, u_2}$
such that
\begin{itemize}
\item[(i)]  \label{eq:cardinality}
$\#\mr_1^{u_1, u_2} > \rho^{8\eta-d}$.
\item[(ii)] If $\|\psi_i-\id\|_{C^1}<\rho$ for $i\in{1,2}$, then
\begin{equation} \label{eq:disjoint}
\left\{ \Pi_{\lambda^{u_1,u_2}}\left(Q^\psi_{u_1^\star 1^\infty,
u_2^\star 1^\infty}(v_1, v_2)\right) : (v_1^\star, v_2^\star)\in
\mr_1^{u_1, u_2}\right\}
\end{equation}
is a pairwise disjoint family, where $Q^\psi_{\omega_1,\omega_2}$ denotes
the rectangle relative to the limit geometry
$\Lambda_{\omega_1,\omega_2}$ (or rather the pair of iterated
function systems defining it).
\end{itemize}

We will later construct a family $\mr^{u_1, u_2}\subset
\mr_1^{u_1,u_2}$ such that $\#\mr^{u_1,u_2}>\rho^{9\eta-\rho}$.
Assuming such a family is given, we define the set of offspring of
$(u_1, u_2)$ to be
\[
V(u_1, u_2) = \{ (u_1 v_1, u_2 v_2) : (v_1^\star, v_2^\star) \in
\mr^{u_1, u_2} \}.
\]

Properties (A) and (C) of $\mathcal{T}$ are clear from the
construction. Property (B) also follows since all $(u_1, u_2)\in
\mathcal{T}_k$ are obtained by going to the $\rho$-decomposition
and then rescaling back to the unit square $k$ times. We will now
consider property (D); along the way we will define the family
$\mr^{u_1,u_2}$ precisely.

Notice that from (\ref{eq:lambda}) and induction we get that for
all $k$ and all $(u_1, u_2)\in \mathcal{T}_k$,
\begin{equation} \label{eq:valuelambda}
\lambda^{u_1, u_2} = \frac{|I(u_2)|}{|I(u_1)|} \lambda.
\end{equation}
Arguing as in the proof of Lemma \ref{lemma:linear}(iii), we get
\[
I^i_{u_i^\star 1^\infty}(v_i)  =  T^i_{u_i} f^i_{u_i} f^i_{v_i}
(T^i_{u_i}f^i_{u_i})^{-1}(I) = T^i_{u_i} f^i_{ u_i v_i}(I),
\]
whence
\[
Q^\psi_{u_1^\star 1^\infty, u_2^\star 1^\infty}(v_1, v_2) =
\psi_1 T^1_{u_1} I^1(u_1 v_1) \times \psi_2 T^2_{u_2}  (I^2(u_2 v_2)).
\]

Since the family in (\ref{eq:disjoint}) is pairwise disjoint, it follows
from (\ref{eq:valuelambda}) that (for fixed $(u_1, u_2)\in
\mathcal{T}_k$) the family
\begin{equation} \label{eq:disjoint2}
\left\{ |I^1(u_1)| \psi_1 T^1_{u_1} I^1(u_1 v_1) + \lambda |I^2(u_2)|
\psi_2 T^2_{u_2} (I^2(u_2 v_2)): (v_1^\star, v_2^\star)\in
\mr_1^{u_1, u_2} \right\}
\end{equation}
is also pairwise disjoint.

For $i\in\{1,2\}$, let $S_i(x) = \mu_i x + \tau_i$ be the
positively-oriented affine map such that $S_i\phi_i$ fixes
$I^i(u_i)$. Since $\|\phi_i-\id\|_{C^1}<\rho/2$, straightforward
calculations and Lemma \ref{lemma:norm} show that
$|\mu_i-1|<\rho/2$ and, restricted to $I^i(u_i)$,
$\|S_i\phi_i-\id\|_{C^1}<\rho$.

Now let
\[
\psi_i= T^i_{u_i} S_i\phi_i (T^i_{u_i})^{-1}.
\]
Notice that $\psi_i\in\mg$ and, by the previous remarks and Lemma
\ref{lemma:norm}, indeed $\psi_i\in\mg(\rho)$. Observe also that
$\{|I^i(u_i)| T^i_{u_i}\}_{i=1,2}$ are translation maps. We deduce
that
\[
|I^i(u_i)| \psi_i T^i_{u_i} = |I^i(u_i)| T^i_{u_i} S_i \phi_i =
\mu_i\phi_i + \tau'_i,
\]
for some $\tau'_i\in\mathbb{R}$. Since affine images of pairwise
disjoint families are still pairwise disjoint, we conclude from
(\ref{eq:disjoint2}) that the following family is pairwise
disjoint as well:
\[
\left\{ \phi_1(I^1(u_1 v_1)) + \frac{\lambda \mu_2}{\mu_1} \phi_2
(I^2(u_2 v_2)): (v_1^\star, v_2^\star)\in \mr_1^{u_1, u_2}
\right\}.
\]
Note that, since $|\mu_i-1|<\rho/2$,
\[
\left|\frac{\mu_2}{\mu_1}-1\right| < \frac{\rho}{1-\rho}.
\]
Hence $|\mu_2\lambda/\mu_1-\lambda|<2A\rho$ whenever $\rho<1/2$,
and it follows from Lemma \ref{lemma:rhoperturb} (or its proof)
that, provided $\rho$ is small enough, there exists a subfamily $\mr^{u_1, u_2}\subset
\mr_1^{u_1,u_2}$ such that
\[
\#\mr^{u_1, u_2} \gtrsim \mr_1^{u_1, u_2} > \rho^{9\eta-d},
\]
and
\begin{equation}  \label{eq:disjoint3}
\left\{\Pi_\lambda Q^\phi(u_1 v_1, u_2 v_2) : (v_1^\star,
v_2^\star)\in \mr^{u_1, u_2} \right\}
\end{equation}
is a pairwise disjoint family. This completes the proof of
Proposition \ref{prop:gugu2}. \qed

\subsection{Conclusion of the proof}
We now complete the proof of Theorem \ref{th:gugu}. First of all
notice that we can assume that $d<1$; if $d\ge 1$ just throw away
some maps in the first i.f.s. (after a suitable iteration) to
obtain a subset of $\Lambda$ of dimension less than, but
arbitrarily close to, 1.

Assume first that $f^i_1$, $i=1,2$, are linear maps, and that
$\log r_1/\log r_2\notin \mathbb{Q}$, where $r_i$ is the
similarity ratio of $f^i_1$. Fix $\varepsilon>0$, and let
\[
S = \left\{ \lambda: \dim_H(\Pi_\lambda(\phi_1(K^1) \times \phi_2(
K^2)))
> d-\e \, \forall \phi_i\in\mg(\delta) \right\}.
\]

By Proposition \ref{prop:gugu2}, if $\delta$ is small enough then
$S$ contains some open set $U$ intersecting both the positive and
negative half-lines. Now let $k,l\in\mathbb{N}$. The cylinders
$K^1(1^{k})$, $K^2(1^{l})$ are, by hypothesis, affine images of
$K^i$ with scaling factors $r_i^k, r_i^l$. This implies that $S$
contains the scaling of $U$ by $\pm r_2^{l}/r_1^{k}$ (the sign
depending on the orientation of the maps $f^i_1$). But $U$ meets
both the positive and negative half-lines, and the set of all such
scaling factors is dense by the irrationality assumption, so we
conclude that $S=\mathbb{R}^*:=\mathbb{R}\backslash\{0\}$.

Next we drop the hypothesis that $f^1_1$, $f^2_1$ are linear; we
still assume that $\log r_1/\log r_2\notin \mathbb{Q}$, where
$r_i$ is the eigenvalue of $f^i_1$ at its fixed point. By Lemma
\ref{lemma:linear}(ii) and the above, there is $\delta>0$ such
that if $\|\phi_i-\id\|_{C^1}<\delta$ for $i\in\{1,2\}$ and $\lambda\in\mathbb{R}^*$, then
\[
\dim_H\left(\Pi_\lambda(\phi_1 L_{1^\infty}(K^1)\times \phi_2
L_{1^\infty}(K_2))\right)
> d-\e.
\]
Let $k$ be so large that $\|T_{1^k}^i
f^i_{1^k}L_{1^\infty}^{-1}\|_{C^1}<\delta$, $i\in\{1,2\}$. Then
\[
\dim_H\left(\Pi_\lambda\left(T^1_{1^k} K^1(1^k) \times T^2_{1^k}
K^2(1^k)\right)\right)
> d -\e\,\textrm{ for all } \lambda\in\mathbb{R}^*.
\]
But since this holds for all $\lambda\neq 0$, $T^i_{1^k}$ is
linear and non-degenerate, and $K^i(1^k)\subset K^i$ for $i\in\{1,2\}$, we obtain
that
\[
\dim_H(\Pi_\lambda(K^1 \times K^2)) > d -\e\,\textrm{ for all }
\lambda\in\mathbb{R}^*.
\]
Since $\e$ was arbitrary, this concludes the proof. \qed

\bigskip

\noindent{\bf Acknowledgements}. I thank C.G.T. de A. Moreira for
explaining the ideas of the proof of Theorem \ref{th:gugu} to me,
and E. J\"{a}rvenp\"{a}\"{a} for useful comments and corrections
on an early version of the article.

This note was written while I was a postdoc at the University of Jyv\"{a}skyl\"{a}. I acknowledge financial support from the Academy of Finland.

%\bibliography{gugu}
%\bibliographystyle{alpha}

\end{document}